\documentclass[11pt]{article}

\usepackage{amssymb}
\usepackage{amsmath}
\usepackage{amsthm}
\usepackage{enumerate}

\newtheoremstyle{plain}%
    {8pt plus2pt minus4pt}%
    {8pt plus2pt minus4pt}%
    {\itshape}%
    {}%
    {\bfseries\scshape}%
    {}%
    {6pt}
    {}%

\newtheoremstyle{remark}%
    {8pt plus2pt minus4pt}%
    {8pt plus2pt minus4pt}%
    {\upshape}
    {}%
    {\bfseries\scshape}%
    {}%
    {6pt}
    {}%

\theoremstyle{plain}
\newtheorem{theorem}{Theorem}[section]
\newtheorem{corollary}[theorem]{Corollary}

\newtheorem{problem}[theorem]{Problem}

\theoremstyle{remark}

\title{Hypergraph Ramsey numbers: tight cycles versus cliques}

\author{
{\Large{Dhruv Mubayi}}\thanks{
\footnotesize {Department of Mathematics, Statistics, and Computer Science, University of Illinois, Chicago, IL 60607, \texttt{Email:mubayi@uic.edu}}. Research partially supported by NSF grants DMS-0969092 and DMS-1300138}
\and
{\Large{Vojtech R\"odl}}\thanks{Department of Mathematics and Computer Science,
Emory University,
Atlanta, GA 30322,
\texttt{Email:rodl@mathcs.emory.edu}. Research partially supported by NSF grants DMS-1102086 and DMS-1301698}}

\date{\today}

\begin{document}
\maketitle		
		
\begin{abstract}
For $s \ge 4$, the 3-uniform tight cycle $C^3_s$ has vertex set corresponding to $s$ distinct points on a circle and edge set given by the $s$  cyclic intervals of three consecutive points.  For fixed $s \ge 4$ and $s \not\equiv 0$ (mod 3)  we prove that there are positive constants $a$ and $b$ with
$$2^{at}<r(C^3_s, K^3_t)<2^{bt^2\log t}.$$  The lower bound is obtained via a probabilistic construction. The upper bound for $s>5$ is proved by using supersaturation and the known upper bound for $r(K_4^{3}, K_t^3)$, while for $s=5$ it follows from a new upper bound for
$r(K_5^{3-}, K_t^3)$ that we develop.
\end{abstract}

\section{Introduction}
A $k$-uniform hypergraph $H$ ($k$-graph for short) with vertex set $V(H)$ is a collection of $k$-element subsets of $V(H)$ (when $k=3$ we will sometimes refer to a 3-graph as a triple system). Write $K_n^k$ for the complete $k$-graph with vertex set of size $n$ (when $k=2$ we omit the superscript for $K_n$ and other well-known graphs). 
 Given $k$ graphs $F, G$, the Ramsey number $r(F,G)$ is the minimum $n$ such that every red/blue coloring of $K_n^k$ results in a monochromatic red copy of $F$ or a monochromatic blue copy of $G$. In this paper we consider hypergraph Ramsey number of cycles versus complete graphs.

Let $C_s$ denote the cycle of length $s$. For fixed $s \ge 3$ the  Ramsey number $r(C_s, K_t)$ has been extensively studied. The case $s=3$ is one of the oldest questions in graph Ramsey theory and it is known that
\begin{equation} \label{r3s} r(C_3,K_t)=\Theta(t^2/\log t).\end{equation}
 The upper bound was proved by Ajtai-Koml\'os-Szemer\'edi~\cite{AKS} and the lower bound was proved by Kim~\cite{K}.  More recently, Bohman-Keevash~\cite{BK2} and independently
Fiz Pontiveros-Griffiths-Morris~\cite{FGM} proved a lower bound that is asymptotically within a factor of 4 of the best upper bound due to Shearer~\cite{Sh}.  The next case $r(C_4, K_t)$ seems substantially more difficult and the best known upper and lower bounds are $O(t^2/\log^2t)$ and $\Omega(t^{3/2}/\log t)$, respectively.
An old open problem of Erd\H os asks whether there is a positive $\epsilon$ for which 
$r(C_4, K_t) = O(t^{2-\epsilon})$.  For larger cycles, the best known bounds can be found in~\cite{BK, SU}. The order of magnitude of $r(C_s, K_t)$ is not known for any fixed $s \ge 4$.

There are several natural ways to extend the definition of cycle to hypergraphs. Two possibilities are to consider loose cycles and tight cycles. The 3-uniform loose triangle $T^3_3$ has 3 edges that have pairwise intersections of size one and have no point in common.  A recent result   due to Kostochka-Mubayi-Verstra\"ete~\cite{KMV} is that there are positive constants $c_1, c_2$ with
$$\frac{c_1 t^{3/2}}{(\log t)^{3/4}} < r(T_3^{3}, K_t^{3})< c_2 t^{3/2}.$$
The authors in \cite{KMV} conjectured that $r(T_3^{3}, K_t^{3})=o(t^{3/2})$.  For longer loose cycles, the gap between the upper and lower bound is much greater.

 For $s \ge k$, the tight cycle $C_s^k$ is the $k$-graph
with vertex set $Z_s$ (integers modulo $s$) and edge set 
$$\{\{i, i+1, \ldots, i+k-1\}: i \in Z_s\}.$$  We can view the vertex set of $C_s^k$ as $s$ points on a circle and the edge set as the $s$ circular subintervals each containing $k$ consecutive vertices. 
In this note we investigate the hypergraph Ramsey number $r(C_s^3, K_t^3)$ for fixed $s\ge 5$ and large $t$. 

When $s \equiv 0$ (mod 3) the tight cycle $C_s^3$ is 3-partite, and in this case it is trivial to observe that $r(C_s^3, K_t^3)$ grows like a polynomial in $t$. Determining the growth rate of this polynomial appears to be a very difficult problem. We focus on the case  $s \not\equiv 0$ (mod 3). In this regime we show that the Ramsey number is exponential in $t$. 

\begin{theorem} \label{main}
Fix $s \ge 5$ and $s \not\equiv 0$ (mod 3). There are positive constants $a$ and $b$ such that
$$2^{at}<r(C^3_s, K^3_t)<2^{bt^2\log t}.$$ 
\end{theorem}

Our upper bounds for $s \ge 7$ are proved by observing that the well-known supersaturation phenomenon in extremal hypergraph theory can be naturally extended to hypergraph Ramsey problems. This observation may be of independent interest. Our upper bound for $s=5$ is by proving an upper bound for $r(K_5^{3-}, K_t^3)$ that matches (apart from constants in the exponent) the best bound for $r(K_4^3, K_t^3)$.  We do this by showing  that the method of Conlon-Fox-Sudakov~\cite{CFS} applies to a particular graph in the ordered setting.  

\section{Lower bound construction}
Our proof is based on the following random construction.

Consider the random graph $G=G(n,1/2)$ with vertex set $[n]$.  Form the (random) triple system $H=H(n)$ where $V(H)=[n]$ and for $i<j<k$, we have
$ijk \in H$ iff $ij, ik \in G$ and $jk \not\in G$.

{\bf Claim.}  If $C^3_s \subset H(n)$, then $3|s$. 

{\bf Proof of Claim.} Assume  that $V(C^3_s)=\{v_1, \ldots, v_s\}$, the edge set is $\{v_iv_{i+1}v_{i+2}:i \in [s]\}$, where indices are taken modulo $s$, and $v_1=\min_i v_i$.
We will prove by induction on $i$ that 
$$v_iv_{i+1} \not\in G \hbox{ iff } i\equiv 2 \hbox{ (mod 3)}.$$
  This is true for $1 \le i \le 3$ as $v_1v_2v_3 \in H(n)$ so $v_1v_2 \in G$ and 
$v_2v_3 \not\in G$. As $v_2v_3v_4 \in H(n)$ we must have $v_4< \min\{v_2, v_3\}$ and hence $v_3v_4 \in G$.  Now suppose that $i\equiv 0$ (mod 3), the statement is true for $i-2, i-1$, and $i$ and let us show it for $i+1, i+2, i+3$.  Since $v_{i-1}v_iv_{i+1} \in H(n)$ and $v_{i-1}v_i \not\in G$, we have $v_{i+1}< \min\{v_{i-1}, v_i\}$ and $v_iv_{i+1} \in G$. Together with $v_iv_{i+1}v_{i+2} \in H(n)$ this gives $v_{i+1}v_{i+2} \in G$ and $v_{i+1}<v_{i+2}$.
Next, $v_{i+1}v_{i+2}v_{i+3} \in H(n)$ implies that $v_{i+2}v_{i+3} \not\in G$ and $v_{i+1}< \min\{v_{i+2}, v_{i+3}\}$. Finally, $v_{i+2}v_{i+3}v_{i+4} \in H(n)$ and $v_{i+2}v_{i+3} \not\in G$ implies that $v_{i+3}v_{i+4} \in G$. 
As $v_1v_{s-1}v_s \in H(n)$ and $v_1=\min_i v_i$
 we also have  $v_{s-1}v_{s} \not\in G$. Consequently, $s-1 \equiv 2$ (mod 3) as desired. \qed

A set $T$ of $t$ vertices in $[n]$ is an independent set in $H(n)$ iff no triple of $T$ is an edge of $H(n)$. Let $T$ be such an independent set in $H(n)$ and $S$ be a partial Steiner triple system with vertex set $T$. It is well-known that such $S$ exist with $|S| = (1+o(1))t^2/6$.  Since $T$ is an independent set, all triples of $S$ are absent in $H(n)$.
As no two triples of $S$ share two points, the probability that all triples of $S$ are absent in $H(n)$ is at most 
$(7/8)^{|S|}<(7/8)^{t^2/7}$. So the probability that there exists a $t$-set  $T$ that is independent in $H(n)$ is at most
$${n \choose t} (7/8)^{t^2/7} < 1$$ as long at 
$t>c\log n$ for some constant $c>0$. We conclude that there exists an $H=H(n)$ with $\alpha(H) < c\log n$, and consequently for fixed $s \not\equiv 0$ (mod 3), there is a constant $c>0$ such that  
$$R(C_s^3, K_t^3) > 2^{ct}.$$

\section{Supersaturation}
Given a hypergraph $H$ and vertex $v \in V(H)$, we say that $w \in V(H)$ is a {\it clone} of $v$ if no edge contains both $v$ and $w$ and for every $e \in H$,
 $v \in e$ if and only if $(e\cup\{w\})\setminus v \in H$.
Given a triple system $F$ and a vertex $v$ in $F$, let $F(v)$ be the triple system obtained from $F$ by replacing $v$ with two clones $v_1,v_2$.
\begin{theorem}  \label{super}
Let $F$ be a triple system with $f$ vertices and $v \in V(F)$. Then  
$$r(F(v), K_t^3)< (r(F, K_t^3))^{f}.$$
\end{theorem}
\proof Let $m=r(F, K_t^3)$ and $n=m^{f}$. Given a red/blue coloring of $K_n^3$, we will find a red copy of $F$ or a blue copy of $K_t^3$. Assume  that there is no blue copy of $K_t^3$. For every $m$-set of vertices, we either find a blue copy of $K_t^3$ or a red copy of $F$ by definition of $m$. Consequently, we may assume that the number of red copies of $F$ is at least 
${n \choose m}/{n-f \choose m-f}= (n)_f/(m)_f> {n \choose f-1}$ since $n=m^f$. To each red copy 
of $F$ associate the $f-1$ vertices  in this copy that do not play the role of $v$. By the pigeonhole principle, we obtain two red copies of $F$ with the same copy of $F-v$, and this yields a red copy of $F(v)$. \qed

A blowup of $F$ is a a hypergraph obtained by successively applying the cloning operation. 
In particular, if each vertex is cloned $p$ times, then denote the obtained blowup as $F(p)$. By applying Theorem \ref{super} repeatedly, we obtain the following easy corollary. 

\begin{corollary} \label{supercor}
Fix a $k$-graph $F$ and an integer $p \ge 2$.
There exists  $c=c(F,p)$ such that 
$$r(F(p), K_t^3) < (r(F, K_t^3))^c.$$
\end{corollary}

Conlon-Fox-Sudakov~\cite{CFS} proved that 
$r(K_4^3, K_t^3)= 2^{O(t^2\log t)}$ and it is an easy exercise to see that for $s \ge 6$, we have $C_s^3 \subset K_4^3(s)$. In other words, there is a homomorphism from $C_s^3$ to $K_4^3(s)$. Consequently, Corollary \ref{supercor} implies the  upper bound in Theorem \ref{main} for $s \ge 6$.  One can check that there is no homomorphism from $C_5^3$ to $K_4^3$ so we cannot apply the same argument for $C_5^3$. 
Clearly 
$$r(C_5^3, K_t^3) < r(K_5^3, K_t^3) < 2^{O(t^3 \log t)}$$
where the last inequality was proved in \cite{CFS}. In the next section we will improve this bound, and this will yield an alternative proof for all $s \ge 5$. 

\section{$K_5^3$ minus an edge}
In this section we use the arguments from \cite{CFS}  to prove that
\begin{equation} \label{k5-} r(K_5^{3-}, K_t^3) = 2^{O(t^2\log t)}\end{equation}
where $K_5^{3-}$ is the triple system obtained from $K_5^3$ by deleting one edge. 
Since $C_5^3 \subset K_5^{3-}$ we immediately obtain the upper bound in Theorem \ref{main} for $C_5^3$.

Given an ordering of the vertices $v_1< v_2<\cdots$ and an ordered graph $F_<$, define the {\it vertex online ordered Ramsey number of $F_<$} as follows: Consider the
following game, played by two players, builder and painter: at step $i+1$ a new vertex $v_{i+1}$ is revealed;
then, for every existing vertex $v_j$, 
$j = 1, \ldots, i$, builder decides, in order, whether to draw the edge
$v_jv_{i+1}$; if he does expose such an edge, painter has to color it either red or blue immediately. The
vertex on-line ordered Ramsey number 
$r(F_<, K_t)$ is then defined as the minimum number of edges that builder
has to draw in order to force painter to create a red ordered $F_<$ (the ordering is defined by the order in which the vertices are exposed) or a blue $K_t$ (since all edges are present in $K_t$, the ordering of its vertices  does not matter). Vertex online ordered Ramsey numbers were studied recently in~\cite{CFLS}.

Let $K_4^-$ be the ordered graph with four vertices $v_1<v_2<v_3<v_4$ and all edges except $v_3v_4$. In other words, it is a copy of $K_4$ with the  ``last" edge deleted.
Our main result is the following.

\begin{theorem} \label{k4-}
In the vertex on-line ordered Ramsey game, builder has a strategy which ensures a red $K_4^-$ or a
blue $K_{t-1}$ using at most
$\ell=O(t^2)$ vertices, $a=O(t^2)$ red edges, and $m=O(t^3)$  
 total edges.
\end{theorem}
Using this result, we immediately obtain the following Theorem which was  proved in the unordered setting in~(\cite{CFS}, Theorem 2.1). The ordered setting we consider here requires essentially no changes to the proof, so we do not repeat the argument.

\begin{theorem} \label{game} Suppose in the vertex on-line ordered Ramsey game that builder has a strategy which ensures
a red $K_4^-$ or a blue $K_{t-1}$ using at most 
$\ell$ vertices, $a$ red edges, and in total $m$ edges. Then, for any
$0 < \alpha < 1/2$, the Ramsey number 
$r(K_5^{3-}, K_t^3)$ satisfies
$$r(K_5^{3-}, K_t^3)= O(\ell \alpha^{-a}(1-\alpha)^{a-m}).$$
\end{theorem}
Using $\alpha=1/t$ and plugging in the values for $\ell,a,m$ from Theorem \ref{k4-} into Theorem \ref{game} yields (\ref{k5-}).

\bigskip

{\bf Proof of Theorem \ref{k4-}.} The idea of the proof is almost identical to Lemma 2.2 in~\cite{CFS}, except that the conclusion one obtains is a bit stronger. In order to apply it to our situation, we need to make a few adjustments. The proof is based on defining builder's strategy which ensures a red graph $K_4^-$ or a blue $K_{t-1}$ in the required number of steps. During the game, builder chooses the vertices one by one and exposes some edges one by one (in the way we will describe below) between the last chosen, say $v_{i+1}$, and the previously chosen vertices $v_1, \ldots, v_i$. The decision of which edge to expose next will depend on the colors of previously chosen edges.  Now we describe the strategy in more detail.

Builder will assign strings consisting of $R$'s and 
$B$'s to each chosen vertex and these strings are expanded during the game in the following way. Assume that the vertices $v_1, \ldots, v_i$ have been chosen and edges between these vertices have been exposed.  After choosing $v_{i+1}$, the first edge exposed is (always) $v_1v_{i+1}$.  Depending on whether painter colors $v_1v_{i+1}$ red or blue the first digit in the string assigned to $v_{i+1}$ will be $R$ or $B$ respectively.  For $j$, $2 \le j \le i$, the only reason why the edge $v_jv_{i+1}$ is exposed is that the current labels of $v_j$ and $v_{i+1}$ are identical.  Let $X_1X_2\ldots X_p$ with $X_p \in \{R, B\}$ for $1 \le q \le p$ be that label.  In order to ``resolve the tie" between the labels of $v_j$ and $v_{i+1}$ builder relables $v_{i+1}$ by a new label which is either $X_1X_2\ldots X_pR$ or $X_1X_2\ldots X_p B$ depending on whether painter colored the edge $v_jv_{i+1}$ red or blue.  Builder stops choosing the vertices and exposing the edges once painter has completed a red graph $K_4^-$ or blue $K_{t-1}$.  It remains to show that during the game at most $\ell=O(t^2)$ vertices, at most $O(t^2)$ red edges, and at most $O(t^3)$ total edges were exposed.

For a vertex $v=v_j$ and color $X \in \{R, B\}$, write $N_X(v)$ for the set of vertices $y=v_i$, $i>j$ that are exposed by builder after $v$ such that painter has colored edge $vy$ with $X$.  In other words, it is the ``forward" neighborhood of $v$ in color $X$.  An important observation we will use in the rest of the proof is the following 

{\bf Fact.} If $y$ is the first exposed vertex in $N_X(v)$, then all edges of the form  $yz$ with $z \in N_X(v)\setminus \{y\}$ are drawn by builder (and consequently colored by painter).

It is convenient to use the following two claims in the proof.

{\bf Claim 1.} Suppose that  
$|N_R(v)|\ge 2s$ for some $v$ and $s\in [t-1]$. Then there is a red $K_4^-$ within $N_R(v) \cup \{v\}$, or a 
blue $K_{s}$ within $N_R(v)$.

{\bf Proof of Claim.} 
Let $w_{a}$ be the first chosen vertex in $W_1=N_R(v)$.
If $w_{a}$ has two red neighbors $x < y$ in $W_1$,
 then $v<w_{a}<x<y$ is a red ordered $K_4^-$ as desired. Consequently, $w_{a}$ has at most one vertex  in $W_1$ (call it $w_{a'}$ if it exists) joined to it in red. Let $W_2=W_1\setminus\{w_a, w_{a'}\}$ so   $|W_2|\ge 2(s-1)$. Note that the Fact implies that all vertices in $W_2$ (including the first) are joined to $w_{a}$ (the first vertex of $W_1$),  by a blue edge. Continuing this process, we either find a red $K_4^-$ or obtain a sequence of proper subsets $W_1\supset W_2 \supset \cdots \supset W_{s}$ with $|W_j|\ge 2(s-j+1)$ for each $j \in [s]$. The first vertex of each $W_j$ for $1 \le j \le s$  forms a set that induces a blue $K_{s}$ within $N_R(v)$ as desired.\qed

{\bf Claim 2.}  Suppose that $1 \le i \le t-2$ and $|N_B(v)|\ge 2{i+1\choose 2}$. Then there is a red $K_4^-$ within $N_B(v)$, or a 
blue $K_{i+1}$ within $\{v\} \cup N_B(v)$.

{\bf Proof of Claim.}  Let us proceed by induction on $i$. The base case $i=1$ is trivial as 
$|N_B(v)|\ge 2>1$ and $v$ with the first exposed  vertex in $N_B(v)$ forms a blue $K_2$.  Now suppose that the result holds for $i-1$ and we wish to prove it for $i \ge 2$. Let $w_a$ be the first chosen vertex in $N_B(v)$. By the Fact, all edges of the form $w_a w_{a'}$ for $w_{a'} \in N_B(v)\setminus\{w_a\}$ have been drawn by builder and hence have been colored, so $w_a$ has at least $2{i+1 \choose 2}-1=2{i \choose 2}+2i-1$ neighbors in $N_B(v)\setminus\{w_a\}$.  By the pigeonhole principle, either at least $2i$ of these neighbors are joined to $w_a$ with a red edge or at least $2{i\choose 2}$ of these neighbors are joined to $w_a$ with a blue edge. In the former case, Claim 1 yields either a red $K_4^-$ or a blue $K_i$ within $N_B(v)$. Together with $v$ this gives a blue $K_{i+1}$ as desired. 
In the latter case, by induction we find a blue $K_{i}$ within $N_B(v)$ and together with $v$ this yields a blue $K_{i+1}$ as desired.
 \qed

Consider the online ordered game where the number of exposed vertices is $2{t \choose 2}=2{t-1 \choose 2}+2(t-1)$.  The number of vertices adjacent to the first vertex $v_1$ is $2{t-1 \choose 2}+2(t-1)$ so at least $2(t-1)$ of them are joined to $v_1$ with a red edge or at least $2{t-1 \choose 2}$ of them are joined to $v_1$ with a blue edge.  In the former case, apply Claim 1 with $s=t-1$ and in the latter case, apply Claim 2 with $i=t-2$.  The claims imply that we find a red $K_4^-$ or a blue $K_{t-1}$. It follows from Claim 2 applied with $i=t-2$ that the number of vertices chosen is at most
 $\ell\le 2{t-1 \choose 2}+1$. Each vertex $v$ is joined to at most three  vertices before it in red, else we have a red $K_4$ (which contains a red $K_4^-$), so the number of red edges is at most $3\ell= O(t^2)$.
Similarly, each vertex is joined to at most $t-2$
vertices before it in blue, else we have a blue $K_{t-1}$. Hence $m\le (t+1)\ell = O(t^3)$ as desired.  
\qed

\section{Concluding remarks}

$\bullet$ The arguments we presented that showed 
$r(K_5^{3-}, K_t^3)=2^{O(t^2\log t)}$ can easily be extended to show that for each fixed $s \ge 4$ we have
\begin{equation} \label{generals}r(K_s^{3-}, K_t^3)=2^{O(t^{s-3}\log t)}.\end{equation}
Note that the  best known bound for fixed cliques of size $s\ge 4$ is $r(K_{s}^{3}, K_t^3)=2^{O(t^{s-2}\log t)}$~\cite{CFS}. Consequently, 
 (\ref{generals}) means that deleting one edge from $K_s^3$ results in a substantially better bound.
  This leads us to ask the following question.

\begin{problem}  Let $s \ge 5$ be fixed.  Is 
$$\log r(K_{s-1}^{3}, K_t^3) = o(\log (r(K_s^{3-}, K_t^3))?$$
\end{problem}

$\bullet$  As $s$ gets large, the tight cycle $C_s^3$ becomes sparser, so one might expect that the Ramsey number with $K_t^3$ (for fixed $s$) decreases as a function of $t$. We were not able to show this and pose the following problem.

\begin{problem} 
Is the following true?
For each fixed $s \ge 4$ there exists  $\epsilon_s$ such that $\epsilon_s \rightarrow 0$ as $s \rightarrow \infty$ and
$$r(C_s^3, K_t^3) < 2^{t^{1+\epsilon_s}}$$
for all sufficiently large $t$.
\end{problem}

$\bullet$ To address the problem of determining $r(C_s^k, K_t^k)$ for fixed $s>k>3$, we can generalize the construction that we used for $k=3$.  
Consider the random $(k-1)$-graph $G^{k-1}=G^{k-1}(n,1/2)$ with vertex set $[n]$.  Form the (random) $k$-graph $H=H(n)$ where $V(H)=[n]$ and for $i_1<i_2<\cdots < i_k$, we have
$e=\{i_1,i_2,\ldots, i_k\} \in H$ iff 
$\{i_2, \ldots, i_k\} \not\in G^{k-1}$ and all other $(k-1)$-subsets of $e$ are in $G^{k-1}$.
Then the arguments we gave for $k=3$ show that 
$C_s^k \in H$ iff $s \equiv 0$ (mod $k$), and the independence number of $H$ is $O((\log n)^{1/(k-2)})$.  Consequently, for fixed $s > k$ with
$s \not\equiv 0$ (mod $k$)
there is a positive constant $c_k$ such that  
$$r(C_s^k, K_t^k) > 2^{c_kt^{k-2}}$$
for all $t$.
The best upper bound we are able to obtain (for fixed $s>k$ and all $t$) is from $r(K_s^k, K_t^k)$. Define the tower function by $t_1(x)=x$  and $t_{i+1}(x)=2^{t_i(x)}$ for $i \ge 1$. Then $r(K_s^k, K_t^k)< t_{k-1}(t^{d_s})$ for all fixed $s > k$, where $d_s>0$ depends only on $s$. 
Consequently, we have
$$2^{c_kt^{k-2}} < r(C_s^k, K_t^k) <r(K_s^k, K_t^k)< t_{k-1}(t^{d_s}).$$
 Closing the  gap above seems to be a very interesting open problem.  In fact, even the first case $s=k+1$, which is probably the easiest as far as lower bounds are concerned, is not known to grow like a tower of height $k-1$.  For example, it is not known whether $r(K_5^4, K_t^4)=r(C_5^4, K_t^4)$ grows as a double exponential (tower of height 3) in $t$. Indeed, finding the minimum $s$ such that $r(K_s^k, K_t^k)$ grows like a tower of height $k-1$ is an open problem raised by Conlon-Fox-Sudakov~\cite{CFS2}.  They proved that this minimum $s$ is at most $\lceil 5k/2\rceil-3$ and conjectured that this can be improved to $s=k+1$. 
Recently, a group at a  workshop on Graph Ramsey Theory at the American Institute of Mathematics observed that a lower bound for $r(K_{k+1}^k, K_t^k)$ can be obtained from a construction that avoids tight paths of appropriate sizes in the ordered setting (they were concerned with the case $t=k+2$ but the same observation is valid for large $t$, which is what we are concerned with here). 
   Results of Shapira-Moshkovitz~\cite{MS} and Milans-Stolee-West~\cite{MSW} provide such a construction and consequently imply that for fixed  $k\ge 3$ and  all $t$
 sufficiently large,
$$ r(C_{k+1}^k, K_t^k) =r(K_{k+1}^k, K_t^k) > t_{k-2}(bt)$$
for some positive constant $b$.  Indeed, the above results of~\cite{MS, MSW} imply that a lower bound for $r(K_{k+1}^k, K_t^k)$ is one more than the size of the poset $J_k$ defined inductively as follows: $J_1$ consists of two chains, one of size 1 and the other of size $t-k$, and for $i\ge 1$, the elements of $J_{i+1}$ are the ideals (down sets) of $J_i$ and order is defined by containment. It is an easy exercise to see that $J_3$ contains an antichain of size $\Omega(t)$ and consequently that $|J_4|=2^{\Omega(t)}$ and then $|J_{\ell}| =t_{\ell-2}(\Omega(t))$ for all $\ell>0$. Similar results can be proved using the ideas developed much earlier in~\cite{DLR}. This suggests that the growth rate of
$r(C_s^k, K_t^k)$ for all fixed $s>k$ is of tower type.

\bigskip

{\bf Acknowledgment.}  The first author thanks David Conlon for an informative short conversation about $r(K_5^{3-}, K_t^3)$.  Part of the research for this paper was done at the workshop Graph Ramsey Theory sponsored by the American Institute of Mathematics in January 2015. The authors are grateful to AIM for hosting the workshop.

\end{document}